\documentclass[reqno,11pt]{amsart}
\usepackage{amsmath}
\usepackage{times}
\usepackage[tcidvi]{graphicx}
\usepackage{color}

\newtheorem{theorem}{Theorem}

\newtheorem{corollary}[theorem]{Corollary}

\newtheorem{definition}[theorem]{Definition}

\newtheorem{lemma}[theorem]{Lemma}

\begin{document}
\title{Relaxation time is monotone in temperature in the
mean-field Ising model}
\author{Vladislav Kargin}
\thanks{Department of Mathematics, Stanford University, CA 94305;
kargin@stanford.edu}
\date{February 2011}
\maketitle

\begin{center}
\textbf{Abstract}
\end{center}

\begin{quotation}
In this note we consider the Glauber dynamics for the mean-field Ising model, when
all couplings are equal and the external field is uniform. It is proved that the relaxation time of the dynamics is monotonically decreasing in temperature.
\end{quotation}

Let $G=\left( V,E\right) $ be a connected graph with $n$ vertices and let $%
\mathcal{S}$ be the set of all assignments of numbers $+1$ or $-1$ to
vertices in $V$. It is convenient to write an element of $\mathcal{S}$ as a
vector $\sigma $ with coordinates $\sigma _{x}=\pm 1,$ where $x\in V.$

The Gibbs measure $\pi $ of the Ising model is a probability measure on $%
\mathcal{S}$: 
\begin{equation*}
\pi (\sigma )=\frac{1}{Z}\exp \left\{ \sum_{x,y\in V}J_{xy}\sigma _{x}\sigma
_{y}+\sum_{x\in V}H_{x}\sigma _{x}\right\} ,
\end{equation*}%
where $Z$ is a normalization factor, $J_{xy}$ are real non-negative numbers
(``couplings''), and $H_{x}$ are real numbers (``external field'')$.$ It is
assumed that $J_{xy}=0$ if $x$ is not connected to $y$ by an edge of the
graph.

We use notation $\left\langle f\right\rangle $ to denote the average of
function $f$ with respect to the Gibbs measure:%
\begin{equation*}
\left\langle f\right\rangle =E_{\pi }f:=\sum_{\sigma \in \mathcal{S}}f\left(
\sigma \right) \pi \left( \sigma \right) .
\end{equation*}

The Glauber dynamics is a reversible Markov chain on $\mathcal{S}$ such that
the Gibbs measure $\pi $ is stationary. Specifically, the transitions
probabilities are as follows. If assignments $\sigma $ and $\sigma ^{\prime
} $ differ on more than one vertex, then $P(\sigma \rightarrow \sigma
^{\prime })=0.$ If they differ on vertex $x$, then 
\begin{equation}
P(\sigma \rightarrow \sigma ^{\prime })=\frac{1}{n}\frac{1}{1+\exp \left(
-2\sigma _{x}^{\prime }\left( S_{x}+H_{x}\right) \right) },
\label{transition_probabilities}
\end{equation}%
where%
\begin{equation*}
S_{x}=\sum_{y\sim x}J_{xy}\sigma _{y}.
\end{equation*}%
Finally, if $\sigma ^{\prime }=\sigma ,$ then 
\begin{equation*}
P(\sigma \rightarrow \sigma ^{\prime })=1-\sum_{x\in V}P(\sigma \rightarrow 
\widehat{\sigma }^{x}),
\end{equation*}%
where $\widehat{\sigma }^{x}$ denote the assignment obtained from $\sigma $
by changing the assignment at vertex $x$.

Let $\mathcal{L}^{2}(\mathcal{S})$ be the linear space of all functions on $%
\mathcal{S}$ with the scalar product 
\begin{equation*}
\left\langle f,g\right\rangle _{\pi }:=\left\langle fg\right\rangle
=\sum_{\sigma \in \mathcal{S}}f(\sigma )g(\sigma )\pi (\sigma ),
\end{equation*}%
where $\pi $ is the Gibbs measure.

Since $P$ is a reversible chain, there is a basis $\{f_{\alpha }\}$ such
that 
\begin{equation*}
Pf_{\alpha }=\lambda _{\alpha }f_{\alpha },
\end{equation*}%
where $1=\lambda _{1}>\lambda _{2}\geq \ldots \geq \lambda _{n},$ and 
\begin{equation*}
\left\langle f_{\alpha },f_{\beta }\right\rangle _{\pi }=\delta _{\alpha
\beta }.
\end{equation*}

The speed of convergence to the equilibrium is governed to a large extent by
the spectral gap $g:=1-\lambda _{2}$. It is conjectured that for every
connected graph $G$ and every family of non-negative couplings $J_{xy}$, 
\begin{equation*}
\frac{\partial \lambda _{2}}{\partial J_{xy}}\geq 0.
\end{equation*}%
See, for example, Question 2 on p. 299 in Levin et al. (2009). This
conjecture has been verified analytically only in the case of the $n$-cycle
with arbitrary couplings (Nacu, 2003). In this note we verify
this conjecture for the case of the mean-field model, in which all the
couplings $J_{xy}$ are the same and the external field is uniform. In this
case formula (\ref{transition_probabilities}) becomes

\begin{equation}
P(\sigma \rightarrow \sigma ^{\prime })=\frac{1}{n}\frac{1}{1+\exp \left(
-2\sigma _{x}^{\prime }\left( J\sum_{y\sim x}\sigma _{y}+H\right) \right) },
\end{equation}%
where $x$ is the only vertex at which $\sigma $ and $\sigma ^{\prime }$ are
different.

In this case, we prove the following result.

\begin{theorem}
\label{theorem_main} Let $G$ be a complete graph on $n$-vertices, let $%
J_{xy}=J>0$ for all $x,y\in V,$ and let $H_{x}=H$ for all $x\in V.$ Let $%
\lambda _{2}$ be the second-largest eigenvalue of the Glauber-Ising model.
Then, it is true that $\lambda _{2}$ is increasing in $J,$%
\begin{equation*}
\frac{\partial \lambda _{2}}{\partial J}\geq 0.
\end{equation*}
\end{theorem}

The relaxation time of the Glauber dynamics is defined as $t_{rel}=\left(
1-\lambda _{2}\right) ^{-1},$ and the temperature $T$ is a parameter
proportional to $J^{-1}.$ Hence, Theorem \ref{theorem_main} has the
following corollary:

\begin{corollary}
For the mean-field model, the relaxation time $t_{tel}$ is decreasing in $T,$%
\begin{equation*}
\frac{\partial t_{rel}}{\partial T}\leq 0.
\end{equation*}
\end{corollary}

The Glauber dynamics on the complete graph was studied as early as in Griffiths et al. (1966), where it was shown that the relaxation time
is exponentially growing in the number of vertices provided that the
temperature is below a critical threshold. Recently, this model has been
investigated in Levin et al. (2010) and in Ding et al. (2009) from the point of view of the theory of finite
Markov chains. There it was shown that the Diaconis cutoff phenomenon (Diaconis (1996)) holds in this model when the temperature is above the
threshold. In addition, these papers investigated the convergence to
equilibrium near the critical temperature and in the slow-convergence regime.

\bigskip

\textbf{Proof of Theorem \ref{theorem_main}}:

\begin{lemma}
\label{lemma_reduction}Let $M$ be the transition matrix of a reversible
Markov chain with stationary distribution $\pi .$ Let $M$ depend on
parameter $J$, and let $\lambda $ be an eigenvalue of $M$ with eigenvector $%
f $ such that $\left\langle f,f\right\rangle _{\pi }=1.$ Then, 
\begin{equation*}
\frac{\partial \lambda }{\partial J}=\left\langle f,\frac{\partial M}{%
\partial J}f\right\rangle _{\pi }.
\end{equation*}
\end{lemma}

\textbf{Proof of Lemma \ref{lemma_reduction}}: We have 
\begin{equation*}
\lambda =\sum_{i,j}f_{i}\pi _{i}M_{ij}f_{j}.
\end{equation*}%
Hence,%
\begin{eqnarray*}
\frac{\partial \lambda }{\partial J} &=&\sum_{i,j}\left( \frac{\partial f_{i}%
}{\partial J}\pi _{i}M_{ij}f_{j}+f_{i}\frac{\partial \pi _{i}}{\partial J}%
M_{ij}f_{j}+f_{i}\pi _{i}\frac{\partial M_{ij}}{\partial J}f_{j}+f_{i}\pi
_{i}M_{ij}\frac{\partial f_{j}}{\partial J}\right) \\
&=&\lambda \frac{\partial }{\partial J}\left( \sum_{i}f_{i}\pi
_{i}f_{i}\right) +\sum_{i,j}f_{i}\pi _{i}\frac{\partial M_{ij}}{\partial J}%
f_{j} \\
&=&\left\langle f,\frac{\partial M}{\partial J}f\right\rangle _{\pi }.
\end{eqnarray*}

QED.

Hence, it remains to prove that $\left\langle f,\frac{\partial P}{\partial J}%
f\right\rangle _{\pi }\geq 0.$ Instead of proving this inequality for the
original Glauber chain, we will define a reduced chain that has the same
second eigenvalue $\lambda _{2}\left( J\right) $ and prove the corresponding
inequality for this new chain. (This chain is called the magnetization chain in %
Levin et al. (2010) and Ding et al. (2009).)

In order to define this new chain, note that every permutation of vertices
induces a linear transformation on $\mathcal{L}^{2}\left( S\right) .$
Because of the symmetry of the mean-field model, the original Glauber chain
has an invariant subspace $L$ that consists of the functions in $\mathcal{L}%
^{2}\left( S\right) $ that are invariant relative to these transformations.
The new transition matrix $\widetilde{P}$ is defined as the restriction of
the original matrix $P$ to this invariant subspace. In more detail, let $%
f\in L$ and let $f_{k}$ be the value of $f$ on configurations with $k$ spins 
$+1$ and $n-k$ spins $-1.$ We will write $f$ as a vector $\left(
f_{0},f_{1},\ldots ,f_{n}\right) .$ This is essentially a choice of a basis
in $L$. Then the transition matrix $\widetilde{P}$ with respect to this
basis is tridiagonal with the entries 
\begin{equation*}
\widetilde{P}_{k,k+1}=\frac{n-k}{n}\frac{1}{1+e^{(n-2k-1)2J-2H}},
\end{equation*}%
where $0\leq k\leq n-1,$ 
\begin{equation*}
\widetilde{P}_{k,k-1}=\frac{k}{n}\frac{1}{1+e^{-(n-2k+1)2J+2H}},
\end{equation*}%
where $1\leq k\leq n,$ and 
\begin{equation*}
\widetilde{P}_{kk}=1-\widetilde{P}_{k,k-1}-\widetilde{P}_{k,k+1},
\end{equation*}%
where $0\leq k\leq n$ and by convention $\widetilde{P}_{0,-1}=\widetilde{P}%
_{n,n+1}=0.$

\begin{definition}
An eigenvector $f=\left\{ f_{k}\right\} _{k=0}^{n}$ of matrix $\widetilde{P}$
is called\emph{\ increasing} if $f_{k+1}\geq f_{k}$ for every $k.$ It is
called \emph{strictly increasing} if it is increasing and $f_{k+1}>f_{k}$
for at least one $k.$
\end{definition}

\begin{lemma}
Matrix $\widetilde{P}$ has the same second-largest eigenvalue $\lambda _{2}$
as $P,$ and this eigenvalue has a striclty increasing right eigenvector.
\end{lemma}

This fact was shown in Ding et al. (2009), in the statement and
proof of Proposition 3.9.

\begin{lemma}
The second largest eigenvalue $\lambda _{2}$ of matrix $\widetilde{P}$ has a
unique increasing eigenvector modulo a multiplication by a scalar.
\end{lemma}

\textbf{Proof:} Let $f=\left\{ f_{k}\right\} _{k=0}^{n}$ be an increasing
eigenvector. We will use the following fact from the proof of Proposition
3.9. in Ding et al. (2009): if $f_{k-1}=f_{k},$ then $%
f_{k-1}=f_{k}=0.$ From this fact it follows that if $f_{k-1}=f_{k}$ and $%
f_{k+1}\neq 0,$ then 
\begin{equation*}
\left( Pf\right) _{k}=\widetilde{P}_{k,k-1}f_{k-1}+\widetilde{P}_{kk}f_{k}+%
\widetilde{P}_{k,k+1}f_{k+1}=\widetilde{P}_{k,k+1}f_{k+1}.
\end{equation*}%
Since $\left( Pf\right) _{k}=\lambda _{2}f_{k}=0,$ hence $\widetilde{P}%
_{k,k+1}=0.$ This contradicts the definition of $\widetilde{P}_{k,k+1}.$
Hence $f_{k}<f_{k+1}$ for every $k.$

Now, let $f$ and $g$ be two increasing eigenvectors corresponding to $%
\lambda _{2}.$ By what we just proved, $f_{k}<f_{k+1}$ and $g_{k}<g_{k+1}$
for every $k.$ Let 
\begin{equation*}
r=\min_{k}\frac{f_{k+1}-f_{k}}{g_{k+1}-g_{k}}.
\end{equation*}
Then $h=f-rg$ is either a zero vector or an increasing eigenvector of $%
\lambda _{2}$ such that $h_{k}=h_{k+1}$ for some $k.$ The latter is
impossible and we showed that modulo a multiplication by a scalar there exists
only one increasing eigenvector of $\lambda _{2}.$ QED.

The symmetry of the model implies that $g_{k}:=\left\{ -f_{n-k}\right\}
_{k=0}^{n}$ is another strictly increasing eigenvector of $\widetilde{P}$
with eigenvalue $\lambda _{2}$. By the previous lemma $g_{k}=f_{k},$ which
means that $f_{k}=-f_{n-k}.$ Since the eigenvector is increasing this
implies that $f_{k}\leq 0$ for $k\leq n/2$ and $f_{k}\geq 0$ for $k\geq n/2.$

Let us define the following quantities: 
\begin{equation*}
s_{k}=\frac{k(n-2k+1)}{n}\frac{1}{1+\cosh \left[ (n-2k+1)2J-2H\right] },
\end{equation*}%
where $0\leq k\leq n.$ Note that $s_{k}\geq 0$ for $k\leq \left( n+1\right)
/2$ and $s_{k}\leq 0$ for $k\geq \left( n+1\right) /2.$

The matrix $\partial \widetilde{P}/\partial J$ is tridiagonal with entries 
\begin{eqnarray*}
\frac{\partial }{\partial J}\widetilde{P}_{k,k+1} &=&s_{n-k}, \\
\frac{\partial }{\partial J}\widetilde{P}_{k,k-1} &=&s_{k}\text{, and} \\
\frac{\partial }{\partial J}\widetilde{P}_{k,k} &=&-s_{k}-s_{n-k},
\end{eqnarray*}%
where $k$ changes between $0$ and $n.$

Let $\widetilde{P}^{\prime }$ denote $\partial \widetilde{P}/\partial J$.
Then we can write: 
\begin{eqnarray*}
f_{0}\left( \widetilde{P}^{\prime }f\right) _{0} &=&f_{0}s_{n}\left(
f_{1}-f_{0}\right) , \\
f_{k}\left( \widetilde{P}^{\prime }f\right) _{k} &=&f_{k}\left[ -s_{k}\left(
f_{k}-f_{k-1}\right) +s_{n-k}\left( f_{k+1}-f_{k}\right) \right] ,\text{ if }%
1\leq k\leq n-1, \\
f_{n}\left( \widetilde{P}^{\prime }f\right) _{n} &=&f_{n}\left[ -s_{n}\left(
f_{n}-f_{n-1}\right) \right] .
\end{eqnarray*}

Since the eigenvector $f$ is increasing, hence all differences $%
f_{k}-f_{k-1} $ are non-negative. Moreover, if $k\leq \left( n-1\right) /2,$
then $f_{k}\leq 0$, $s_{k}\geq 0,$ and $s_{n-k}\leq 0,$ which implies that $%
f_{k}\left( \widetilde{P}^{\prime }f\right) _{k}\geq 0.$ Similarly, $k\geq
\left( n+1\right) /2$ implies that $f_{k}\left( \widetilde{P}^{\prime
}f\right) _{k}\geq 0.$ The only remaining case is when $n$ is even and $%
k=n/2.$ However, in this case $f_{k}=0$ and therefore $f_{k}\left( 
\widetilde{P}^{\prime }f\right) _{k}=0.$ It follows that 
\begin{equation*}
\left\langle f,\frac{\partial \widetilde{P}}{\partial J}f\right\rangle _{\pi
}=\sum_{k=0}^{n}\pi _{k}f_{k}\left( \widetilde{P}^{\prime }f\right) _{k}\geq
0.
\end{equation*}%
This completes the proof of Theorem \ref{theorem_main}.

\end{document}